\documentclass[11pt, a4paper,twoside]{article}

\usepackage{amsmath,amssymb,amsthm,amscd}
\usepackage[]{fontenc}
\usepackage{xy}
\usepackage{enumerate}

\xyoption{all} \numberwithin{equation}{section}

\begin{document}

\title{Slope inequalities for fibred surfaces via GIT}
\author{Lidia Stoppino\footnote{Work partially supported by: 1) PRIN 2003: \textit{Spazi di moduli e teoria di Lie}; 2) GNSAGA. 2000 Mathematics Subject Classification: Primary 14H10, Secundary 14J29, 14D06} }
\date{}
%\footnote{email: stoppino@mat.uniroma3.it}}

\newcommand{\Q}{\mathbb{Q}}
\newcommand{\Z}{\mathbb{Z}}
\newcommand{\R}{\mathcal{R}}
\newcommand{\pr}{\mathbb P}
\newcommand{\sym}{\mbox{\upshape{Sym}}}
\newcommand{\of}{\omega_f}
\newcommand{\rk}{\mbox{\upshape{rank}}}
\newcommand{\av}{``}
\newcommand{\og}{\omega_\beta}
\newcommand{\oa}{\omega_\alpha}
\newcommand{\op}{\omega_{f'}} 
\newcommand{\g}{\gamma}
\newcommand{\y}{{Y_b}}
\newcommand{\el}{\mathcal L}
\newcommand{\oo}{\mathcal O} 
\newcommand{\surjarrow}{\rightarrow\!\!\!\!\rightarrow}
\newcommand{\rd}{\mbox{\upshape{red.deg}}}
\newcommand{\kap}{\mathbb{C}}
\newcommand{\pil}{\pi_{\Lambda}}
\newcommand{\cliff}{\mbox{\upshape{Cliff}}}
\newcommand{\pic}{\mbox{\upshape{Pic}}}
\newcommand{\fidi}{\varphi_{|D|}}
\newcommand{\gon}{\mbox{\upshape{gon}}}
\newcommand{\ci}{c_1}
\newcommand{\ch}{\mbox{\upshape ch}}
\newcommand{\td}{\mbox{\upshape td}}
\newcommand{\kapp}{\mbox{\itshape {\textbf k}}}
\newcommand{\Si}{\Sigma}
\newcommand{\qb}{{q_f}}

\newtheorem{teo}{Theorem}[section]
\newtheorem{prop}[teo]{Proposition}
\newtheorem{lem}[teo]{Lemma}
\newtheorem{rem}[teo]{\it Remark}
\newtheorem{defi}[teo]{\it Definition}
\newtheorem{ex}[teo]{\it Example}
\newtheorem{ass}[teo]{Assumption}
\newtheorem{cor}[teo]{Corollary}

\renewcommand{\theequation}{\arabic{section}.\arabic{equation}}

\pagestyle{myheadings}
\markboth{\small{Lidia Stoppino}}{\small{\textit{Slope inequalities via GIT}}}

\maketitle

\vspace{-0.5cm}

\begin{abstract}
In this paper we present a generalisation of a theorem due to Cornalba and Harris, which is an application of Geometric Invariant Theory to the study of invariants of fibrations. 
In particular, our generalisation makes it possible to treat the problem of bounding the invariants of general fibred surfaces. 
As a first application, we give a new proof of the slope inequality and of a bound for the invariants associated to double cover fibrations. 
\end{abstract}

%\tableofcontents

\section*{Introduction}
In 
\cite{X}
and 
\cite{C-H},
Xiao and Cornalba-Harris 
developed two methods that can be applied to the problem of bounding the invariants  of  fibred varieties. 
Given a complex variety $X$ fibred over a curve,
the starting point of both methods is a line bundle $L$ on $X$.
However, while Xiao's method uses techniques of vector bundle stability, the one of Cornalba-Harris  exploits 
Geometric Invariant Theory (GIT).
In the same papers, the three authors treat the case of surfaces fibred over a curve, 
proving a fundamental inequality on the invariants: the so-called slope inequality (cf. section \ref{bounds}). 
However,  Cornalba and Harris prove the inequality only for   semistable fibrations 
(i.e. fibred surfaces such that any fibre is a semistable curve in the sense of Deligne and Mumford). 
In fact, their method applies only to semistable non-hyperelliptic fibrations, 
and the semistable hyperelliptic case is obtained by an ad hoc argument.

A result of Tan \cite{T1} made apparent that  the general case of the slope inequality  
cannot be reduced to the semistable case (see Remark \ref{seminotenough}).

The starting point of this work is the question whether or not it is possible to treat the non semistable  
and the hyperelliptic cases using the ideas of Cornalba-Harris. 
The answer is affirmative; while the extension to non-semistable fibrations is straightforward, 
in order to treat the  hyperelliptic case it is necessary to modify substantially the method.
This led us to develop a generalisation of the method which is, in our opinion, interesting on its own.

Until now, the method of Xiao have  been almost the only one used to find lower bounds on the slope of fibred 
surfaces (although  the very nice argument introduced by Moriwaki in \cite{Mor} should also be taken into account, 
as remarked also in \cite{A-K}).
It has been further developed and applied by Ohno,  Konno, Barja, Zucconi, and others.

Thanks to the generalisation presented here,  the  Cornalba-Harris method qualifies as  a valid alternative. 
In this paper, besides proving  the slope inequality, we can show in  a new and  direct way that
the fibred surfaces realising the equality are all hyperelliptic. 
Furthermore, we prove an inequality holding between the invariants of double cover fibrations.
%(Definition  \ref{dcf}).

In  \cite{B-S} both the generalised Cornalba-Harris method and the one of Xiao have been applied, 
obtaining a new bound on the slope of non-Albanese fibrations.
Besides, the first method seems to have promising applications to the case of fibrations of higher dimensional varieties 
 \cite{trigonalebarjastoppino}, where the one of Xiao tends to be  technically 
hard.\\
\indent The importance of these kind of results is double. 
On the one hand, they are fundamental tools in the study of the geography of complex  surfaces
(for example, Pardini's recent proof of the Severi inequality in \cite{Parda}, makes an essential use of the slope inequality).
On the other hand, the bounds on  the slope have an application to the positivity of divisors on  the moduli 
space of stable curves of genus $g$
(for instance, in \cite{GKM},  the slope inequality is a key ingredient for attaching a conjecture on the
nef cone).

We now give a brief account of the method of Cornalba-Harris and of its generalisation.
The idea of the method is the following.
Let $f\colon X\rightarrow T$ be a flat proper morphism of complex varieties 
with a line bundle $L$ on $X$ whose restrictions to the general fibres of $f$ give embeddings in projective spaces.
Suppose that the Hilbert points of these embeddings are semistable in the sense of GIT. 
Then the key-point of the method is to translate the  semistability assumption into the existence of a line bundle   
on the base $T$, together with a non-vanishing section of it.
This produces in particular an element in the effective cone of the base $T$. 
When $T$ is a curve, the consequence is a non-trivial inequality holding between  the degrees of rational classes 
of divisors on it.\\
\indent The main point of the generalisation  is to drop the assumption that the line bundle $L$ gives an 
\emph{embedding} on the general fibres, and to consider {\em arbitrary rational maps}.
In order to do this, we need to introduce a suitable generalisation of Hilbert (semi)stability for a 
variety with a  map in a projective space (Definitions \ref{GHS} and \ref{Hs}). 
This generalisation sounds unexpected because,
as GIT is mainly used to construct moduli spaces, 
GIT stability is usually defined for line bundles whose associated morphisms encode 
all the information about the variety, as in the case of the classical Hilbert points.
We prove that, assuming this generalised semistability, the argument of  Cornalba-Harris, 
with some modifications, can be pushed through, and still gives as a consequence an 
effective divisor on the base $T$ (Theorem \ref{CH}).\\
\indent The paper is organised as follows.
In the first section we prove the main theorem (Theorem \ref{CH}); under some more restrictive 
assumptions, we can  derive explicit inequalities on the rational classes of divisors on the base (Corollary \ref{quelo}).
In section \ref{bounds}, we give the proof of the slope inequality, and of the fact that the 
 fibrations with minimal slope are all hyperelliptic (Proposition \ref{caratt}).
We treat in section \ref{DC} the case of surfaces having an involution on a genus $\g$ fibration 
(double cover fibrations, see Definition \ref{dcf}), proving an inequality on the invariants. \\
{\it Acknowledgements.} It is a pleasure to thank Maurizio Cornalba, who has introduced me to this subject, taught me how to deal with it, and patiently helped me through the preparation of this paper.
I am also indebted to  Miguel \'Angel Barja for helpful conversations, to  Rita Pardini for useful comments on a preliminary version of the paper, and to both of them for their encouragement.  
I thank Alessandro Ghigi for his help on  stability issues.

%%%%%%%%%%%%%%%%%%%%%%%%%%%%%%%%%%%%%%%%%%%%%%%%%%%%%%%%%%%%%%%%%%%%%%%%%%%%%%%%%%%%%%%%%%%%%%%%%%%%%%%%%%%%%%%%%%%

\section{The Cornalba-Harris method generalised}\label{metodo}

%\subsection{Stability of maps to $\pr^s$: generalised Hilbert stability}

Let $G$ be a reductive complex algebraic group and $V$ a finite dimensional representation of $G$. 
An element $v\in V$ is said to be GIT \emph{semistable} if the closure of its orbit does not contain $0$, and GIT \emph{stable} if its stabiliser is finite and its orbit closed.
Recall that a necessary and sufficient condition for the semistability of $v\in V$
is the existence of a $G$-invariant non-constant homogeneous 
polynomial $f\in \mbox{Sym}(V^{\vee})$ such that $f(v)\not =0$.

Let $X$ be a variety (an integral separated scheme of finite type over $\kap$), 
with a linear system $V\subseteq H^0(X,L)$, 
for some line bundle  $L$ on $X$.
Fix $h\geq 1$ and call $G_h$ the image of the natural homomorphism
\begin{equation}\label{omo}
\sym^hV\stackrel{\varphi_h}{-\!\!\!\longrightarrow}H^0(X, L^h).
\end{equation}
Set $N_h=\dim G_h$ and take exterior powers
$$
\bigwedge^{N_h}\sym^hV\stackrel{\wedge ^{N_h}\varphi_h}{-\!\!\!\longrightarrow}
\bigwedge^{N_h}G_h=\det G_h.
$$
If we identify $\det G_h$ with $\kap$, the homomorphism $\wedge ^{N_h}\varphi_h$ can be seen as a liner functional on $\wedge^{N_h}\sym^hV$. Changing the isomorphism, it gets multiplied by a non-zero element of $\kap$. Hence, we can see $\wedge ^{N_h}\varphi_h$ as a well-defined element of $\pr (\wedge^{N_h}\sym^hV^{\vee})$.

\begin{defi}\label{GHS}
With the above notations, we call 
$\wedge ^{N_h}\varphi_h\in \pr(\wedge^{N_h} \mbox{Sym}^hV^{\vee})$, 
the  generalised $h$-th Hilbert point 
associated to the couple  $(X, V)$. 
\end{defi}
If $V$ induces an embedding, then for $h\gg 0$ the homomorphism $\varphi_h$ 
is surjective  and  
it is the classical $h$-th Hilbert point associated to $\psi$. 

Let $\dim V=s+1$ and consider the  standard representation  $SL(s+1,\kap)\rightarrow SL(V)$; 
we get  an induced natural action of $SL(s+1,\kap)$ on  $\pr(\wedge^N \mbox{Sym}^hV^{\vee})$, 
and we can introduce the associated notion of GIT (semi)stability:
we say that the $h$-th generalised Hilbert point of the couple  $(X, V)$
is {\em semistable (resp. stable) if it is GIT semistable (resp. stable) with respect to the natural 
$SL(s+1, \kap)$-action.}

\begin{rem}\label{stabimmagine}{\upshape
Let   $(X, V)$ be as above.
Consider the factorisation of the induced map through the image,
$$X-\rightarrow \overline X \stackrel{j}{\hookrightarrow}\pr^s.$$
Set $\overline L=j^*(\oo_{\pr^s}(1))$ and let $\overline V\subseteq H^0(\overline X,\overline L)$
 be the linear systems associated to $j$.
The homomorphism (\ref{omo}) factors as follows:
$$\sym^hV\cong\sym^h\overline V\stackrel{\overline \varphi_h}{\longrightarrow} H^0(\overline X,\overline{L}^h)
\hookrightarrow H^0(X, L^h),$$
where the homomorphism $\overline \varphi_h$ is the $h$-th Hilbert point of the 
embedding $j$; 
notice that, by Serre's vanishing theorem, this homomorphism is onto  
(and, in particular, $G_h=H^0(X,\overline{L}^h)$)  for large enough $h$.
The generalised $h$-th Hilbert point of $(X, V)$ is therefore naturally identified 
with the $h$-th Hilbert point of $(\overline X,\overline V)$, and  the generalised Hilbert stability 
of $(X,V)$ coincides with the classical Hilbert stability 
of the embedding $j$. }
\end{rem}
\begin{defi}\label{Hs}
We say that $(X,V)$ is generalised Hilbert stable (resp. semistable) if its generalised $h$-th Hilbert point is stable 
(resp. semistable) for infinitely many integers $h>0$.
\end{defi}
In the case of embeddings in projective space, this notion coincides with the classical 
Hilbert stability introduced in \cite{Mum}.

\subsection{The theorem}
We will use the following well known fact about vector bundles and representations.
\begin{rem}\label{vbr}
{\upshape
Let $T$ be a  projective variety.
Consider a vector bundle $E$ of rank $r$ on $T$ and a complex holomorphic 
representation 
$$
GL(r,\kap)\stackrel{\rho }{\longrightarrow} GL(V).
$$
Composing the transition functions of $E$ with $\rho$, we can construct a new vector bundle, which we call $E_\rho$.
Hence,
if $\{ g_{\alpha,\beta}\}$ is a system of transition functions for $E$ with respect to an open
cover $\{\mathcal U_{\alpha, \beta}\}$ of $T$, then a system of transition functions for $E_\rho$ with respect 
to the same cover is $\{\rho(g_{\alpha, \beta})\}$. 
Clearly $E_{\rho}$ has typical fibre $V$. 
%and structure group $\rho (GL(r,\kap))$. 

For instance, if we consider as $\rho$ the representation corresponding to symmetric, tensor and exterior power, the vector bundle $E_\rho$ becomes respectively 
$\sym^n E$, $\otimes ^nE$ and $\wedge^nE$.}
\end{rem}

%%%%%%%%%%%%%%%%%%%%%%%%%%%%%%%%%%%%%%%%%%%%%%%
\noindent We are now ready to state the theorem. Notation: given a sheaf $\mathcal F$ 
over a variety $T$, 
we call $\mathcal F\otimes \kapp(t)$ the fibre of $\mathcal F$ over the point $t\in T$.

\begin{teo}\label{CH}
Let $f\colon X\rightarrow T$ be a flat morphism from a variety $X$ to a variety $T$. Let 
$t$ be a general point of $T$, $X_t$ the fibre of $f$ at $t$. 
Let $L$ be a line bundle 
on $X$ and $\mathcal F$ a locally free subsheaf of $f_*L$ of rank $r$.
Suppose that  for some integer $h>0$ the $h$-th  generalised Hilbert point associated to the linear system
$\mathcal F\otimes \kapp(t)\subseteq H^0(X_t,L_{|X_t})$ 
is semistable.
Let $\mathcal G_h\subseteq f_*L^h$ be a locally free subsheaf that contains the image
of the morphism 
$$
\mbox{\upshape{Sym}}^h\mathcal F\longrightarrow f_*L^h,
$$
and coincides with it at $t$.
Set $N_h=\rk \,\mathcal G_h$.
Let $\mathcal L_h$ be the  line bundle
$$
\mathcal L_h=\det(\mathcal G_h)^{r}\otimes (\det \mathcal F)^{-h\,N_h}.
$$

Then there is a positive integer $m$, depending only on $h$, $\rk \mathcal F$ and $N_h$, such that $(\mathcal L_h)^m$ is effective.
\end{teo}

\begin{proof}
In what follows, $t$ is a general point of $T$.
Set $F:=\mathcal F\otimes \kapp(t)$, $G_h:=\mathcal G_h\otimes \kapp (t)$.
Consider the morphism
$
\mbox{Sym}^h\mathcal F\stackrel{\gamma_h}{\longrightarrow}\mathcal G_h.
$
Its fibre at $t$, 
$\overline \gamma_h\colon \mbox{Sym}^h  F \longrightarrow G_h,$
is surjective by assumption.
Its maximal exterior power is the generalised Hilbert point associated to  $(X_t, F)$.
Therefore there exists by assumption a  homogeneous $SL(F)$-invariant polynomial (of degree, say, $d$)
$P\in \mbox{Sym}^d(\bigwedge^{N_h}\mbox{Sym}^h F)$
such that 
\begin{equation}\label{P}
\sym^d\bigwedge^{N_h}\overline \gamma_h(P)\not= 0 \,\,\,\mbox{ in }\,\,\,(\det G_h)^d.
\end{equation}
We may assume (simply taking a power of $P$ if necessary) 
that the degree of $P$ is  $mr$, where $m$ is an integer
depending only on $h$, $r$ and $N_h$.
Fixing an isomorphism $F \cong \kap^r$, $P$ corresponds to an element $$\widetilde P\in \sym^{mr}\bigwedge^{N_h}\sym^n \kap^r.$$
If we change the isomorphism, as  $P$ is  invariant by the action of $SL(F)$, we obtain $\widetilde P$ multiplied by a non-zero element of $\kap$. 
Hence, the line $W$ generated by $\widetilde P$ in $\mbox{Sym}^{mr}(\wedge^{N_h}\mbox{Sym}^h \kap^r)$,  is well defined and  invariant under the action of $GL(r,\kap)$.
This produces naturally a line bundle on $T$ with an injective morphism into $(\det \mathcal G_h)^{mr}$, as we verify at once, using the language of representations.

Let  $\rho$ be the $N_h$-th exterior power of the $h$-th symmetric 
power of the standard representation, 
$$
\rho \colon GL(r,\kap)\rightarrow GL(\wedge^{N_h}\mbox{Sym}^h \kap^r)).
$$ 
Using the notations of  Remark \ref{vbr}, the vector bundle  $\mathcal F_\rho$ is $\wedge^{N_h}\mbox{Sym}^h\mathcal F$.
Let $$\sigma\colon GL(r,\kap)\longrightarrow GL(W)$$ 
be the representation obtained by restriction from $\sym^{mr}\rho$. 
Thus, there is an inclusion of vector bundles
$
\mathcal F_\sigma\hookrightarrow \sym^{mr}\mathcal F_\rho.
$
Composing this inclusion  with $\mbox{Sym}^{mr}\wedge^{N_h}\gamma_h$, we obtain  a homomorphism 
$\mathcal F_\sigma \rightarrow (\det \mathcal G_h)^{mr}$, whose fibre at $t$ is the following composition
$$
W\hookrightarrow \sym^{mr}\bigwedge^N\sym^n( F)\longrightarrow (\det  G_h)^{mr},
$$
which is a non-zero homomorphism by construction because of property (\ref{P}) (it is, roughly speaking, the evaluation of  $\gamma_h$ on $P$).
It remains to understand explicitly $\mathcal F_\sigma$.
Given an element $M\in GL(r,\kap)$, if we write 
$M=(\det M)^{1/r}U$, where $U\in SL(r,\kap)$, the action of $M$ on $P$ is the following:
$$
\sigma (M) P=\mbox{Sym}^{mr}\rho((\det M)^{1/r}U) P=\det \mu (M)^{hNm}\mbox{Sym}^{mr}\rho(U)P=\det \mu (M)^{hNm}P.
$$
It follows that in our case $\mathcal F_\sigma$ is the line bundle $(\det \mathcal F)^{hNm}$, and the proof  is  concluded.
\end{proof}

%%%%%%%%

%\subsection{Applications to the effective cone of $T$}\label{applications}

In all the applications of the Cornalba-Harris method that have been made so far, including ours, 
the condition of stability is 
satisfied not for a {\em fixed $h$}, but for {\em $h$ large enough}: 
more precisely Hilbert stability is satisfied (see Definition \ref{Hs}).

Moreover, it is often the case that the choice of $\mathcal F\subseteq f_*L$ and of $\mathcal G_h$ is such that 
the first rational Chern class  $\ci(\mathcal L_h)\in A^1(T)_\Q$ is a polynomial in $h$
of the form
\begin{equation}\label{poli}
\ci(\mathcal L_h)=\alpha_d h^d+\ldots +\alpha_1 h + \alpha_0,\,\,\,\,\,\,\,\alpha_i\in A^1(T)_\Q.
\end{equation}
Theorem \ref{CH} assures that for infinitely many positive integers $h$ there exists an integer 
$m$ such that the line bundle $\el_h^m$ is effective, 
hence the class $\ci(\mathcal L_h)\in  A^1(T)_\Q$ is effective.
In this situation, 
we can therefore conclude that  the leading coefficient $\alpha_{d}$  
is the limit in $A^1(T)_\Q$ of effective divisors.\footnote{Note that, 
although we are speaking of limits, we don't need to pass to real coefficients, 
because in fact both $\alpha_d$ and the members of the succession converging to 
it given by (\ref{poli}) belong to $A^1(T)_\Q$.}

We can make explicit computations and simplifications under additional assumptions (this corollary should be compared to  
the original Theorem (1.1) of \cite{C-H}). 

\begin{cor}\label{quelo}
With the notations of Theorem \ref{CH}, suppose  that $\mathcal F$ induces a Hilbert semistable map on the general fibres.
Suppose moreover that

(1) $f$ is proper, $T$ is irreducible of dimension $k$ and $X$ is of pure dimension $k+d$;

(2) for $t\in T$ general, the fibre $\mathcal F\otimes\kapp(t)$ induces an embedding of $X_t$;

(3) the higher direct images $R^if_*L^h$ vanish for $i>0$, $h \gg 0$ (this happens for instance if the fibre of $\mathcal F$
induces an ample linear system on {\em any} fibre of $f$).

Then, the class 
$$\mathcal E(L,\mathcal F):=rf_*(\ci(L)^{d+1}\cap [X])-
(d+1)\ci (\mathcal F)\cap f_*(\ci (L)^d\cap [X])$$
  is contained in the closure of the effective cone of $A_{k-1}(T)_\Q$.
\end{cor}
\begin{proof}
By the second  assumption, for general $t$, the homomorphism
$$
\sym ^h\mathcal F\otimes \kapp(t)\longrightarrow H^0(X_t, L^h_{|X_t})
$$ 
is surjective for large enough $h$. 
Hence, we can choose $\mathcal G_h=f_*L^h$ in Theorem \ref{CH}.
Therefore, 
$$
\ci(\mathcal L_h)= r\,\ci(f_*L^h)-h\,\rk f_*L^h\,\ci(\mathcal F).
$$
The first assumption enables us to use  the Riemann-Roch Theorem for singular varieties 
(\cite{fulton}, Corollary 18.3.1) and obtain the formula 
\begin{equation}\label{GRR}
\ch \left(f_! L^h\cap \td(\oo_T)\right) =f_*\left(\ch (L^h)\cap \td (\oo_X)\right).
\end{equation}
%where 
%$$\ch(\mathcal E)=\rk \mathcal E\, [X]+\ci (\mathcal E)+\frac{1}{2}(\ci(\mathcal E)^2-2 c_2(\mathcal E))+\ldots\in A_*(X)_\Q$$and  
%$$\td (\mathcal E)=[X]+\frac{1}{2}\ci(\mathcal E)+\frac{1}{12}(\ci(\mathcal E)^2+ c_2(\mathcal E))+\ldots\in A_*(X)_\Q$$ 
%are respectively the Chern character and the Todd class of a sheaf $\mathcal E$ on $X$, and 
%$$f_! \mathcal E=\sum_i(-1)^i R^if_*\mathcal E.$$
%Moreover, recall that, as $L$ is a line bundle,
%$$\ch (L^h)=\sum_{i=0}^{k+d}\frac{\ci(L^h)^i}{i!}=\sum_{i=0}^{k+d}\frac{h^i\ci(L)^i}{i!}.$$
Recalling that, for any variety $Y$, $\td (\oo_Y)=[Y]+ \mbox{ terms of dimension }< \dim Y$,
and using standard intersection-theoretical computations, we obtain that
\begin{equation}\label{cia}
\begin{array}{cl}
\ci(\mathcal L_h)\cap [T] = &  \frac{h^{d+1}}{(d+1)!}\mathcal E(L,\mathcal F)+\\
& \\
& +\sum_{i=1}^{d}(-1)^{i+1}\left( r\ci (R^if_*L^h) \cap [T] -h \,\rk(R^if_*L^h)\, \ci(\mathcal F)\cap  [T]\right)+\\
& \\
& + \,O(h^d).
\end{array}
\end{equation}
Hence, equation (\ref{cia}), together with the remarks made above, implies the statement.
\end{proof}

%%%%%%%%%%%%%%%%%%%%%%%%%%%%%%%
%%%%BOUND SULLA SLOPE!!!!!!%%%%%%
%%%%%%%%%%%%%%%%%%%%%%%%%%%%%%%%%

\section{Bounds on the slope of fibred surfaces}\label{bounds}

A \emph{fibred surface}, is the datum of  a surjective 
morphism $f$ with connected fibres from a smooth projective surface $X$ to a smooth complete curve $B$. 
Throughout this section, we shall use the term ``fibration'' as a synonimous of fibred surface.
The genus $g$ of the general fibre is called genus of the fibration.
%
%Define a $-1$-curve (respectively a $-2$-curve) to be a nonsingular rational curve 
%$C\subset X$ with self-intersection $-1 $
%(respectively $-2$).
%
We call a fibration \emph{relatively minimal} if the fibres contain no $-1$-curves.
A fibration is said to be \emph{semistable}  if all the fibres are semistable curves in the sense of Deligne and Mumford
(i.e. if it is relatively minimal with nodal fibres).
From any fibred surface  $f\colon X\rightarrow B$, by  contracting all the $-1$-curves in the fibres, 
we obtain  an induced fibration on $B$, called the relatively minimal model of $f$, which is unique if $g\geq 1$.
We say that a fibration is \emph{locally trivial} if it is a holomorphic fibre bundle.

As usual, the \emph{relative canonical
sheaf} of a fibred surface $f\colon X\rightarrow B$ is
the line bundle $\of =\omega_X\otimes (f^*\omega_B)^{-1}$; and let $K_f$ denote any associated divisor.
From now on we will consider {\em relatively minimal fibrations of genus $g\geq 2$.}
Two basic invariants for such a fibration are the following.
$$K_f^2=K^2_X-8(g-1)(g(B)-1);$$
$$\chi_f=\chi (\oo_X)-\chi (\oo_B)\chi (\oo_F)=\chi(\oo_X)-(g-1)(g(B)-1).$$ 
Using Riemann-Roch and Leray's spectral sequence, one sees that $\chi_f= \deg f_*\of$. 
It is well known that both these invariants are non-negative. 
Moreover, $\chi_f=0$ if and only if $f$ is locally trivial.
Assuming that the fibration is not locally trivial, we can consider the ratio
$$
s(f):= \frac{K_f^2}{\chi_f},
$$
which is called the \emph{slope}. 
Of course $s(f)\geq 0$; but a bigger bound holds, given by the following result, 
which we call  \emph{slope inequality}:
\begin{teo}[Xiao, Cornalba-Harris]\label{teoslope}
Let $f\colon X\longrightarrow B$ be a relatively minimal fibration of genus $g\geq 2$.
\begin{equation}\label{slope}
gK_f^2\geq 4(g-1)\chi_f
\end{equation}
\end{teo}
This inequality is sharp, and it is possible to classify the fibrations reaching it, 
which are in particular all  hyperelliptic (Proposition \ref{caratt}).

\begin{rem}\label{seminotenough}
{\upshape
As is well-known, the process of \emph{semistable reduction} 
associates to any fibred surface a semistable one, by means of a ramified base change. 
One might hope that, using semistable reduction, it could be possible to reduce the proof of the 
slope inequality for any fibration to the semistable case.
However, Tan has shown (cf. Theorem A and Theorem B of \cite{T1}) that the behaviour of the slope under base change \emph{cannot be controlled when the base change ramifies over fibres which are not D-M semistable}, which is precisely what happens in the semistable reduction process.
In particular, the inequalities that can be shown to hold for semistable fibrations, do not necessarily extend to arbitrary fibrations.}
\end{rem}

The form of Theorem \ref{CH}  we shall use in  the applications to surfaces is the following.
\begin{cor}\label{princ1}
Let $f\colon X\rightarrow B$ be a fibred surface. 
Let $L$ be a line bundle on $X$ and $\mathcal F$ a coherent\footnote{As the base $B$ is a smooth curve, any coherent subsheaf of a locally free sheaf is locally free.} subsheaf of $f_*L$ of rank $r$ such that for general $b\in B$ the linear system 
$$\mathcal F\otimes \kapp(b)\subseteq H^0(X_b,L_{|X_b})$$
induces a Hilbert semistable map. 
Let $\mathcal G_h$  be a coherent subsheaf of $f_*L^h$ that contains the image
of the morphism 
$\mbox{\upshape{Sym}}^h\mathcal F\longrightarrow f_*L^h,$
and coincides with it at general $b$.
If $N=\rk \mathcal G_h$ is of the form $A h+ O(1)$ and $\deg \mathcal G_h$ of the form
$B h^2+ O(h)$, the following inequality holds:
\begin{equation}\label{princ2}
rB-A\deg \mathcal F\geq 0.
\end{equation}
\end{cor}
\begin{proof}Straightforward from Theorem \ref{CH} and the observations made  after it.
\end{proof}
%%%%%%%%%%%%%%%%%%%%%%%%%%%%%%%%%%%%
%%%%%%%%%%%%%%%%%%%%%%%%%%%%%%%%%%%%%%%%%%%%%%%%%%%%%%%
We now come to the proof of the slope inequality.
\begin{proof} {\it of Theorem \ref{teoslope}}\\
\noindent We want to apply Corollary \ref{princ1} with $L=\of$ and $\mathcal F=f_*L$. 
Let  $X_b$ be a general fibre. 

Observe that the higher direct image  $R^1f_*\of ^h$ vanishes for large enough $h$,
as can be seen for instance using the relative version of Kawamata-Viehweg vanishing theorem (cf. \cite{KMM}, Theorem 1.2.3).
We split the proof in two steps:\\

\noindent{\it (1) Suppose $f$ is non-hyperelliptic}.
The condition of Corollary \ref{princ1} is satisfied, because the canonical embedding of a smooth non-hyperelliptic curve is Hilbert stable, as shown in \cite{Mum}:  indeed (using Mumford's  notations), it is linearly stable, and hence  Chow stable, which in turns implies the generalised Hilbert stability; see also \cite{ACG2}  or \cite{tesi} for a direct proof.
We can compute the terms in inequality (\ref{princ2}) as follows
$$
\begin{array}{l}
\displaystyle \mbox{rank }f_*\of=h^0(X_b,{\of}_{|X_b})=g;\smallskip \\
\displaystyle \rk \, \mathcal G_h=h^0(X_b,{\of}_{|X_b}^h)=(2h-1)(g-1);\smallskip \\
\displaystyle \deg \mathcal G_h=\frac{(hK_f\cdot (h-1)K_f)}{2}+\deg f_*\of= h^2\frac{K_f^2}{2}+O(h).
\end{array}
$$
Hence, inequality (\ref{princ2}) becomes
%$$g K_f^2\geq 4(g-1)\deg f_*\of =4(g-1)\chi_f,$$which is 
exactly the slope inequality.\\

\noindent{\it (2) Suppose  $f$ is hyperelliptic}.
A general hyperelliptic fibred surface is not always a double cover of a 
fibration of genus $0$.
Anyway we show below that for our purposes 
it can be treated as if it were.
We make use of a standard argument (cf.  for instance \cite{A-K}) which can be applied to 
any fibred surface with an involution that restrict to an involution on the general fibres.
% (see [Bar] and [A-K] for similar constructions).

First observe that  the hyperelliptic involution on the general fibres extends to a global involution 
$\iota$ on $X$ (see for instance \cite{Persdouble}).
If $\iota$ has no isolated fixed points then  
$X/\!\!<\!\iota\!>$ is a smooth genus $0$ surface on $B$ and the quotient map is a double cover 
whose ramification divisor is the fixed locus of $\iota$.
Otherwise, we blow up the isolated fixed 
points and obtain a smooth surface $\widetilde{X}$ birational to $X$ whose 
induced involution $\tilde\iota$ has no isolated fixed points. 
Call $Y$ the quotient of $\widetilde X$ by $\tilde \iota$. 
The surface $Y$ has a natural genus $0$  fibration $\alpha$ over $B$, but is not necessarily
relatively minimal. 
We have the following diagram:
\begin{equation}\label{diag}
\xymatrix{
\widetilde X \ar[dr]^\pi \ar[d]^\eta \\
X  \ar@{-->}[r]  \ar[d]^f & Y\ar[dl]^\alpha\\
B}
\end{equation}
Let $R\subset Y$ be the branch divisor of $\pi$.
By the theory of cyclic coverings (cf. \cite{BHPVdV} I.17), 
we can find a line bundle $\mathcal L$ on $Y$ such that 
$\mathcal L^2=\mathcal O_Y(R)$. 
Set $\tilde f=f\circ\eta$. 
Recall that $\omega_{\tilde f}=\eta^*\of\otimes\mathcal O_{\widetilde X}(E)$, where $E$ is the
 union of the exceptional $-1$-curves. 
Let $\epsilon$ be the number of connected components of $E$. 
%(i.e.the number of blow ups of which $\eta$ is made of). 
Consider the exact sequence 
$$
0\rightarrow \eta^*\of^h\rightarrow \omega_{\tilde f}^h\rightarrow \mathcal O_{hE}(hE)\rightarrow 0
$$
and the long exact sequence induced by the pushforward by $\tilde f$:
$$
0\rightarrow f_*\of^h\rightarrow \tilde f_*\omega_{\tilde f}^h\rightarrow \tilde f_*\mathcal O_{hE}(hE)\rightarrow...
$$
$$
...\rightarrow R^1 f_*\of^h\rightarrow R^1 \tilde f_*\omega_{\tilde f}^h\rightarrow R^1\tilde f_*\mathcal O_{hE}(hE)
\rightarrow 0.
$$
Observe that $\deg \tilde f_*\mathcal O_{hE}(hE)=h^0(\mathcal O_{hE}(hE))= 0$, and that
$$\deg R^1\tilde f_*\mathcal O_{hE}(hE)=
h^1(\mathcal O_{hE}(hE))=\epsilon \frac{h^2-h}{2},$$ 
by the Riemann-Roch Theorem for embedded curves.
Therefore
$
\tilde f_*\omega_{\tilde f}^h=f_*\of^h
$ for any $h$, and
$$
\deg R^1\tilde f_*\omega_{\tilde f}^h=\deg R^1 f_*\of^h+\epsilon \frac{h^2-h}{2}=\epsilon \frac{h^2-h}{2}.
$$
Recall that in our situation  $\omega_{\tilde f}=\pi^*(\oa\otimes\mathcal L)$ and $\pi_*\mathcal O_{\widetilde X}=\mathcal O_Y\oplus \mathcal L^{-1}$. 
Therefore we have the following decomposition of $\tilde f_*\omega_{\tilde f}$
\begin{equation}\label{ugua}
\tilde f_*\omega_{\tilde f}=\alpha_*\pi_*\omega_{\tilde f}=\alpha_*(\pi_*\pi^*(\oa\otimes\el))=\alpha_*((\oa\otimes\el)\otimes\pi_*\oo_Y)=
\alpha_*(\oa\otimes\el)\oplus\alpha_*\oa.
\end{equation}
Hence, $ f_*\omega_{f}=\alpha_*(\oa\otimes\el)$, being $\alpha$  a genus $0$ fibration.

The canonical line bundle
$\omega_{X_b}={\of}_{|X_b}$ induces a morphism to $\mathbb P^{g-1}$ that 
factors through a double cover of $\pr^1$ ramified at the Weierstrass points of $X_b$ composed with the 
Veronese embedding of degree $g-1$. 
The morphism
$
\sym^hf_*\of \longrightarrow f_*\of^h
$
has fibre on $b$
$$
\sym^hH^0(X_b, \omega_{X_b})=\sym^hH^0(\pr^1,\oo_{\pr^1}(g-1))\surjarrow H^0(\pr^1,\oo_{\pr^1}(h(g-1)))\subset H^0(X_b, \omega_{X_b}^h).
$$
Observe that the fibre $\alpha_*(\oa\otimes\el)^h\otimes \kapp (b)$ is
 $H^0(\mathbb P^1,\mathcal O_{\mathbb P^1}(h(g-1)))$;  
we hence choose $\alpha_*(\oa\otimes\el)^h$ as the sheaf $\mathcal G_h$
 in Corollary \ref{princ1}.
The semistability assumption is satisfied, because 
the Veronese embedding $\pr^1\hookrightarrow \pr^{g-1}$ has 
semistable Hilbert point, as shown for instance in \cite{K}, cor. 5.3.
For large enough $h$, by the Riemann-Roch Theorem
$$\deg\mathcal G_h=h^2\frac{(K_\alpha+L)^2}{2} +\deg R^1\alpha_*(\oa\otimes \el)^h+O(h),$$
$$\rk\mathcal G_h=h^0(\y,\omega^h_{\y}(hL))=h(g-1)+1.$$
We now estimate the degree of $R^1\alpha_*(\oa\otimes\mathcal L)^h$ for $h\gg0$. 
Observe that  $R^1\tilde f_*\omega_{\tilde f}^h$ is torsion and 
splits into the direct sum 
$$R^1 \tilde f_*\omega_{\tilde f}^h=R^1\alpha_*(\oa\otimes\mathcal L)^h\oplus R^1\alpha_*
(\oa^h\otimes\mathcal L^{h-1}).$$
Now, observe that for large enough $h$ 
%from the sequence 
%$$0\rightarrow \oa^h\otimes\el^{h-1}\rightarrow \oa^h\otimes\el^h\rightarrow (\oa^h\otimes\el^h)_{\vert L}\rightarrow 0$$
%and the long sequence induced by pushforward $\alpha_*$,
%we see that
$$\deg R^1\alpha_*(\oa\otimes\mathcal L)^h= 
\deg R^1\alpha_*(\oa^h\otimes\mathcal L^{h-1})+O(h),$$
hence
$$\deg R^1\alpha_*(\oa\otimes \el)^h=\frac{1}{2} \deg R^1\tilde f_*\omega_{\tilde f}^h+ O(h)=
\epsilon \, \frac{h^2}{4}+ O(h),$$
and inequality (\ref{princ2}) becomes 
$$
\frac{g}{2}\left((K_\alpha+L)^2+\frac{\epsilon}{2} \right)-(g-1)\deg \alpha_*(\oa\otimes \el)\geq 0.
$$
As   $\pi$ is a finite morphism of degree $2$, and $\eta$ is a sequence of $\epsilon$ blow ups,
$$K_f^2-\epsilon= K_{\tilde f}^2=(\pi^*(K_\alpha+L))^2=
2(K_\alpha+L)^2.$$
Remembering  that 
$\deg \alpha_*(\oa\otimes \el)=\deg \tilde f_*\omega_{\tilde f}=\deg f_*\of=\chi_f$,
we obtain the slope inequality.
\end{proof}

The following proposition has been proved by Konno in \cite{konnohyp}, some years later the 
proof of the slope inequality, as a by-product of other inequalities.
Using the approach of Cornalba-Harris, it is a natural consequence of the construction.
This proposition is a generalisation of the first part of Theorem (4.12) of \cite{C-H}.
\begin{prop}\label{caratt}
Let $f\colon X\rightarrow B$ be a relatively minimal non-locally trivial fibred surface of genus $\geq 2$ satisfying 
equality in Theorem \ref{teoslope}. Then $f$ is hyperelliptic.
\end{prop}
\begin{proof}
Suppose by contradiction that it holds 
$g K_f^2=4(g-1)\chi_f$
and that $f$ is non-hyperelliptic.
Going back to the proof of Theorem \ref{CH}, we see that in order to have equality, 
the degree of the degree $2$ coefficient in the polynomial $\ci(\el_h)$ must be $0$.
As $\ci(\el_h)$ is effective, the linear coefficient has to be of non-negative degree.
Computing this class we get 
$$
0\leq -\frac{g}{2}K_f^2+(g-1)\chi_f=-(g-1)\chi_f,
$$
which is strictly negative for $g\geq 2$ and $f$ non-locally trivial. 
Hence, we get the desired contradiction.
\end{proof}
The hyperelliptic fibrations that reach the bound can be classified, and turn out to have 
restrictions on the type of singularities of the special fibres (see \cite{C-H}, Theorem (4.12) for the 
semistable case, and \cite{A-K} sect. 2.2 for the general one).

%%%%%%%%%%%%%%%%%%%%%%%%%%%%%%%%%%%%%%%%%%%%

\section{Bounds for double cover fibrations}\label{DC}

Arguing in a very similar way to what we did for hyperelliptic fibrations, we can  prove a bound for 
the invariants of a more general class of fibred surfaces, double cover fibrations:
\begin{defi}\label{dcf}
A \emph{double cover fibration of type $(g,\g)$} is the data of a genus $g$ fibred surface 
$f\colon X\rightarrow B$ together with a global involution on $X$ that restricts, on the general fibre, to an involution with genus $\g$ quotient. 
\end{defi}
In particular, the  double cover  fibrations of type $(g,0)$ are exactly the hyperelliptic ones.
The slope of double cover fibrations has been studied in \cite{B-Z}, and recently in \cite{CS}. 
We refer to these two articles for a detailed discussion of the situation. 
In \cite{CS} the sharp bound 
\begin{equation}\label{sharpDC}
s(f)\geq 4\frac{g-1}{g-\g}
\end{equation}
is proved, under the assumption $g\geq 4\g+1$. 
For $g<4\g$ the bound is {\em false} in general. 
Proposition \ref{dc-ch} below implies that the bound holds in general for a particular class of
double cover fibrations.
A similar inequality can be  found applying Xiao's method (\cite{BPhD}, Prop.4.10).

Let  $f\colon X\rightarrow B$ be a double cover fibration of type $(g,\g)$ with $\g\geq 1$.
With the same construction made for the hyperelliptic case, we can associate to it a 
genus $\gamma$ fibration $\alpha\colon Y\rightarrow B$, not necessarily relatively minimal, 
obtaining a diagram of the form (\ref{diag}). 
Let us use the same notations of the hyperelliptic case.

\begin{prop}\label{dc-ch}
Let $f\colon X\rightarrow B$ be a double cover fibration of type $(g,\g)$ with $\g\geq 1$ and  $g\geq 2\g+1$.
Let $\alpha\colon Y\rightarrow B$ be the associated fibration of genus $\g$ described above.

Then the following inequality holds:
\begin{equation}\label{aoL}
K_f^2\geq 4\frac{g-1}{g-\g}(\chi_f- \chi_\alpha).
\end{equation}
In particular, any double cover fibration with $g\geq 2\g+1$ and associated genus $\g$ fibration 
isotrivial satisfies the bound (\ref{sharpDC}).
\end{prop}
\begin{proof}
Arguing as in the hyperelliptic case (same notations), we obtain the decomposition
$$f_*\of= \alpha_* (\oa\otimes \el)\oplus \alpha_*\oa,$$
which on a general fibre $X_b$ amounts to 
$$
H^0(X_b,\omega_{X_b})=H^0(\y,\omega_\y( L))\oplus H^0(\y,\omega_{\y}).
$$
where $L$ is the restriction of $\el$ to $\y$.
By Hurwitz' formula $\deg L=g-2\g+1$.
We want to apply Corollary \ref{princ1} of Theorem \ref{CH} to the rank $g-\g$ subsheaf 
$\mathcal F:=\alpha_*(\oa\otimes\el)$ of $f_*\of$.

We split the proof in  two cases.

(1) Suppose that  the restriction of $\mathcal F$ on a general fibre $\y$ does not belong to a $g^1_2$ on $\y$
(this holds in particular if $\alpha$ is non-hyperelliptic or if $g\geq 2\g+2$).
In this case $\mathcal F$  induces on a general fibre $X_b$ a $2:1$ morphism to $\y$ followed by the 
morphism $\psi$ in $\pr^{g-\g-1}$ induced by the line bundle $\omega_\y(L)$. 
We distinguish again two cases. 
(1.a) $\psi$ is an embedding; in this case
it  is linearly stable,  by \cite{Mum}, section 2.15, 
hence, by the same argument made in the non-hyperelliptic case of Theorem \ref{teoslope}, 
it is Hilbert stable. 
We apply Corollary \ref{princ1} taking as 
$\mathcal G_h$ the sheaf $\alpha_*(\oa^h\otimes\el^h)$.
Now, computing $\deg\mathcal G_h$, $\rk\mathcal G_h$, and 
$\deg R^1\alpha_*(\oa\otimes\mathcal L)^h$ for $h\gg 0$, 
as in the hyperelliptic case of Theorem \ref{teoslope}, inequality (\ref{princ2}) becomes 
$$
\frac{g-\g}{2}\left((K_\alpha +L)^2+\frac{\epsilon}{2} \right)-(g-1)\deg \alpha_*(\oa\otimes \el)\geq 0.
$$
Remembering  that 
$$K_f^2-\epsilon= K_{\tilde f}^2 =\pi^*(K_\alpha + L)^2=
2(K_\alpha + L)^2,$$
and that
$\deg \alpha_*(\oa\otimes \el)=\deg \tilde f_*\omega_{\tilde f}-\deg \alpha_*\oa=\deg f_*\of-\deg \alpha_*\oa=\chi_f-\chi_\alpha,$
we obtain the statement.
(1.b) $\psi$ fails to be an embedding if and only if $\deg L=2$. Note that, by assumption, if $C$ is hyperelliptic, $L\not \in g^1_2$.
In this case $\psi$ is a birational morphism, which
%its image $\overline \y\subset \pr^{\g}$  is an irreducible curve 
%either with one node (if $p\not = q$), or with an ordinary cusp (if $p=q$). 
is linearly semistable, and hence, by \cite{Mum} again, its image is Chow semistable. 
Chow semistability does not imply Hilbert semistability, hence we cannot use the Cornalba-Harris method;
 however, we can in this case apply 
a result of Bost (\cite{bost}, Theorem 3.3) that gives as a consequence exactly the same inequality of Corollary \ref{princ1}.

(2) Suppose on the other hand that $\alpha$ is hyperelliptic and that  the morphism induced by $\alpha_*\oa\otimes \el$ 
on a general fibre factors through the hyperelliptic involution of $\y$:
$$
X_b \stackrel{2:1}{\longrightarrow}  \y \stackrel{2:1}{\longrightarrow} \pr^1\stackrel{v}{\hookrightarrow} \pr^{g-\g-1},
$$
where $v$ is the Veronese embedding.
The semistability assumption is satisfied because $v$ is Hilbert semistable (as observed in the hyperelliptic case in
Theorem \ref{teoslope})
With similar computations, we obtain
$$
\deg \mathcal G_h=h^2\frac{K_f^2}{8}+ O(h),
\quad
\rk \mathcal G_h =\g h+ O(1),
$$
and again inequality (\ref{princ2}) gives the desired  bound.
\end{proof}

\addcontentsline{toc}{section}{References}

\bigskip
\noindent Dipartimento di Matematica, Universit\`a di Pavia, 
Via Ferrata 1, 27100 Pavia - ITALY. \\
E-mail: \textsl {lidia.stoppino@unipv.it}.


\begin{thebibliography}{99}

\bibitem{ACG2} E. Arbarello, M. Cornalba, P. A. Griffiths, J. Harris, \emph{Geometry of algebraic curves. Vol.{II}}, in preparation.

\vspace{-0.1cm}

\bibitem{A-K} T. Ashikaga and K. Konno, \emph{Global and local properties of pencils of algebraic curves}, Algebraic geometry 2000, Azumino (Hotaka), Adv. Stud. in Pure Math. {\bf 36}  Math. Soc. Japan, Tokyo, 2002, 1-49. 

%\bibitem{Ara} S. Arakelov, \emph{Families of algebraic curves with fixed degeneracies}, Math. U.S.S.R. Izv. {\bf 5 } (1971),1277-1302 .

\vspace{-0.1cm}

\bibitem{BPhD} M. A. Barja, \emph{On the slope and geography of fibred surfaces and threefolds}, Ph. D. Thesis, Univesity of Barcelona, 1998.

\vspace{-0.1cm}

%\bibitem{Barja3folds} M. A. Barja, \emph{On the slope of fibred threefolds}, Internat. J. Math. {\bf 11} n.4 (2000), 461-491.

%\bibitem{B} M. A. Barja, \emph{On the slope of bielliptic fibrations}, Proc. Amer. Math. Soc. {\bf 129} n.7 (2001), 1899-1906.

\bibitem{B-S} M. A. Barja and L. Stoppino, \emph{Linear stability of projected canonical curves with applications to the slope of fibred surfaces},  J. Math. Soc. Japan., {\bf 60}, No. 1 (2008) 1-22.

\vspace{-0.1cm}

\bibitem{trigonalebarjastoppino} M.A. Barja and L. Stoppino, \emph{A sharp bound for the slope of general trigonal fibrations of even genus}, in preparation.

\vspace{-0.1cm}


\bibitem{B-Z} M. A. Barja and F. Zucconi, \emph{On the slope of fibred surfaces}, Nagoya Math. J. {\bf 164} (2001), 103-131.

\vspace{-0.1cm}


%\bibitem{Bcan} M. A. Barja {\em Numerical bounds of canonical varieties}, Osaka J. Math. {\bf 37 (3)} (2000), 701-718.

\bibitem{BHPVdV} W. P. Barth, K. Hulek, C. A. M. Peters, A. Van de Ven, \emph{Compact complex surfaces}, second edition, 
Springer-Verlag, Berlin Heidelberg, 2004.

\vspace{-0.1cm}

%\bibitem{Bv} A. Beauville, \emph{L'in\'egalit\'e $p_g\geq 2q+4$ pour les surfaces de type g\'en\'eral}, Bull. Soc. Math. Franc. {\bf 110} (1982), 3
%43-346.
\bibitem{bost} J. Bost, {\em Semi-stability and heights of cycles}, Invent. Math., 
{\bf 118}, vol. 2  (1994), 223-253.



\vspace{-0.1cm}

\bibitem{C-H} M. Cornalba, J. Harris, \emph{Divisor classes associated to families of stable varieties,
with applications to the moduli space of curves.} Ann. Sc. Ec. Norm. Sup. {\bf 21 (4)} (1988), 455-475.

\vspace{-0.1cm}

\bibitem{CS} M. Cornalba and L. Stoppino, \emph{A sharp bound for the slope of double cover fibrations}, to appear in the Michigan Mat. J., preprint math.AG/0510144.

\vspace{-0.1cm}

\bibitem{fulton} W. Fulton, \emph{Intersection Theory, second edition}, Springer-Verlag, 1998.

\vspace{-0.1cm}

\bibitem{GKM} A. Gibney, S. Keel and I. Morrison, {\em Towards the ample cone of {$\overline M\sb {g,n}$}}, J. Amer. Math. Soc.
 {\bf15 (2)}, (2002), 273-294 (electronic).

\vspace{-0.1cm}

\bibitem{KMM} Y. Kawamata and K. Matsuda and K. Matsuki, {\em Introduction to the {M}inimal {M}odel {P}roblem}, Adv. Stud. in Pure Math., {\bf 10}, (1987), 283-360.

\vspace{-0.1cm}

\bibitem{K} G. R. Kempf, {\em Instability in invariant theory}, Ann. of Math. {\bf 108}, (1978), 299-316.

%\bibitem{konnocliff} K. Konno, \emph{Clifford index and the slope of fibered surfaces}, J. Algebraic Geom. {\bf 8 (2)} (1999), 207-220.

%\bibitem{konnotrig} K. Konno, \emph{A lower bound of the slope of trigonal fibrations}, Internat. J. Math. {\bf 7 (1)} (1996), 19-27.

\vspace{-0.1cm}

\bibitem{konnohyp} K. Konno, \emph{Non-hyperelliptic fibrations of small genus and certain irregular canonical surfaces}, Ann. Sc. Norm. Sup. Pisa  ser. IV {\bf 20} (1993), 575-595.

\vspace{-0.1cm}

%\bibitem{KoIrr} K. Konno, \emph{On the {I}rregularity of {S}pecial {N}on-{C}anonical {S}urfaces}, Publ. RIMS Kyoto Univ. {\bf 30} (1994), 671-688.

\bibitem{Mor} A. Moriwaki, {\em A sharp slope inequality for general stable fibrations of curves}, J. Reine Angew. Math.,
{\bf 480},  (1996), 177-195.

\vspace{-0.1cm}

\bibitem{Mum} D. Mumford, \emph{Stability of projective varieties}, L'Ens. Math. {\bf 23} (1977), 39-110.

\vspace{-0.1cm}

\bibitem{Parda} R. Pardini, {\em The {S}everi inequality {$K\sp 2\geq 4\chi$} for surfaces of maximal {A}lbanese dimension},
Invent. Math. {\bf 159 (3)} (2005), 669-672.

\vspace{-0.1cm}

\bibitem{Persdouble} U. Persson, {\em Double coverings and surfaces of general type},
Algebraic geometry (Proc. Sympos., Univ. Troms\o, Troms\o,1977), {Lecture Notes in Math.},
{\bf 687},168--195, 1978.

\vspace{-0.1cm}

\bibitem{tesi} L. Stoppino, \emph{Stability of maps to projective spaces, with applications to the slope of fibred surfaces}, Ph.D. Thesis, {U}niversit\` a di  {P}avia, 2005.

\vspace{-0.1cm}

\bibitem{T1} S. Tan, \emph{On the invariants of base changes of pencils of curves I}, Manuscripta Math. {\bf 84} (1994), 225-244. \emph{On the invariants of base changes of pencils of curves II}, Math. Z. {\bf 222} (1996), 655-676.

\vspace{-0.1cm}

\bibitem{X} G. Xiao, \emph{Fibred algebraic surfaces with low slope.} Math. Ann. {\bf 276} (1987), 449-466.

\end{thebibliography}
\end{document}